\documentclass[draft,11pt]{amsart}
\usepackage{amsmath,amssymb,amscd,amsthm}
\newtheorem{formula}{}[section]
\newtheorem{proposition}[formula]{Proposition}
\newtheorem{corollary}[formula]{Corollary}
\newtheorem{lemma}[formula]{Lemma}
\newtheorem{theorem}[formula]{Theorem}
\theoremstyle{definition}
\newtheorem{definition}[formula]{Definition}
\newtheorem{example}[formula]{Example}
\theoremstyle{remark}
\newtheorem*{remark}{Remark}

\newcommand{\D}{\Delta}
\renewcommand{\l}{\lambda}

\newcommand{\C}{\mathbb C}
\newcommand{\R}{\mathbb R}
\newcommand{\Z}{\mathbb Z}
\renewcommand{\k}{\mathbf k}

\newcommand{\bideg}{\mathop{\rm bideg}}
\newcommand{\Tor}{\mathop{\rm Tor}\nolimits}
\newcommand{\tbtk}{\widetilde{B_TK}}
\newcommand{\zk}{\mathcal Z_K}

\newcommand{\ep}{\nolinebreak\hfill$\square$}

\begin{document}

\title[Torus actions and coordinate subspace arrangements]
{Torus actions, equivariant moment-angle complexes, and coordinate
subspace arrangements}
\author{Victor M. Buchstaber}
\author{Taras E. Panov}
\thanks{Partially supported by
the Russian Foundation for Fundamental Research, grant no. 99-01-00090,
and INTAS, grant no. 96-0770.}
\address{Department of Mathematics and Mechanics, Moscow
State University, 119899 Moscow, RUSSIA}
\email{tpanov@mech.math.msu.su \quad buchstab@mech.math.msu.su}

\begin{abstract}
We show that the cohomology algebra of the complement of a coordinate
subspace arrangement in $m$-dimensional complex space is isomorphic to the
cohomology algebra of Stanley--Reisner face ring of a certain simplicial
complex on $m$ vertices.  (The face ring is regarded as a module over the
polynomial ring on $m$ generators.) Then we calculate the latter
cohomology algebra by means of the standard Koszul resolution of
polynomial ring. To prove these facts we construct an equivariant with
respect to the torus action homotopy equivalence between the complement of
a coordinate subspace arrangement and the moment-angle complex defined by
the simplicial complex. The moment-angle complex is a certain subset of
a unit poly-disk in $m$-dimensional complex space invariant with respect to
the action of an $m$-dimensional torus.
This complex is a smooth manifold provided that the simplicial
complex is a simplicial sphere, but otherwise has more complicated
structure. Then we investigate the equivariant topology of the
moment-angle complex and apply the Eilenberg--Moore spectral sequence.
We also relate our results with well known facts in the theory of toric
varieties and symplectic geometry.
\end{abstract}

\maketitle

\section{Introduction}

In this paper we apply the results of our previous paper~\cite{BP2} to
describing the topology of the complement of a complex coordinate subspace
arrangement. A coordinate subspace arrangement $\mathcal A$ is a set of
coordinate subspaces $L$ of a complex space $\C^m$, and its complement is the
set $U(\mathcal A)=\C^m\setminus\bigcup_{L\in\mathcal A}L$.  The complement
$U(\mathcal A)$ decomposes as $U(\mathcal A)=U(\mathcal A')\times(\C^{*})^k$,
were $\mathcal A'$ is a coordinate arrangement in $\C^{m-k}$ that does not
contain any hyperplane. There is a one-to-one correspondence between
coordinate subspace arrangements in $\C^m$ without hyperplanes and simplicial
complexes on $m$ vertices $v_1,\ldots,v_m$: each arrangement $\mathcal A$
defines a simplicial complex $K(\mathcal A)$ and vice versa. Namely let
$|\mathcal A|$ denotes the support $\bigcup_{L\in\mathcal A}L$ of the
coordinate subspace arrangement $\mathcal A$; then a subset
$v_I=\{v_{i_1},\ldots,v_{i_k}\}$ is a $(k-1)$-simplex of $K(\mathcal A)$ if
and only if the $(m-k)$-dimensional coordinate subspace $L_I\subset\C^m$
defined by equations $z_{i_1}=\ldots=z_{i_k}=0$ does not belong to $|\mathcal
A|$.  An arrangement $\mathcal A$ is obviously recovered from its simplicial
complex $K(\mathcal A)$; that is why we write $U(K)$ instead of
$U\bigl(\mathcal A(K)\bigr)$ throughout this paper. (For more information
about relations between arrangements and simplicial complexes see the
beginning of Section~2.)

Subspace arrangements and their complements play a pivotal role in many
constructions of combinatorics, algebraic and symplectic geometry, mechanics
etc., they also arise as configuration spaces of different classical systems.
That is why the topology of complements of arrangements entranced many
mathematicians during the last two decades.  The first important result here
deals with arrangements of hyperplanes (not necessarily coordinate) in
$\C^m$. Arnold~\cite{Ar} and Brieskorn~\cite{Br} shown that the cohomology
algebra of the corresponding complement $U(\mathcal A)$ is isomorphic to the
algebra of differential forms generated by the closed forms $\frac1{2\pi
i}\frac{dF_A}{F_A}$, where $F_{A}$ is a linear form defining the hyperplane
$A$ of the arrangement. Orlik and Solomon~\cite{OS} proved that the
cohomology algebra of the complement of a hyperplane arrangement depends only
on the combinatorics of intersections of hyperplanes and presented
$H^{*}\bigl(U(\mathcal A)\bigr)$ by generators and relations. In general
situation, the Goresky--MacPherson theorem~\cite[Part~III]{GM} expresses the
cohomology groups $H^{i}\bigl(U(\mathcal A)\bigr)$ (without ring structure)
as a sum of homology groups of subcomplexes of a certain simplicial complex.
This complex, called the {\it order} (or {\it flag}) complex, is defined via
the combinatorics of intersections of subspaces of $\mathcal A$. The proof
of this result uses the stratified Morse theory developed
in~\cite{GM}. Another way to handle the cohomology algebra of the complement
of a subspace arrangement was recently presented by De Concini and
Procesi~\cite{dCP}. They proved that the rational cohomology ring of
$U(\mathcal A)$ is also determined by the combinatorics of intersections.
This result was extended by Yuzvinsky in~\cite{Yu}. In the case of coordinate
subspace arrangements the order complex is the barycentric subdivision of a
simplicial complex $\tilde K$, while the summands in the Goresky--MacPherson
formula are homology groups of links of simplices of $\tilde K$. The complex
$\tilde K$ has the same vertex set $v_1,\ldots,v_m$ as our simplicial complex
$K$ and is ``dual" to the latter in the following sense: a set
$v_I=\{v_{i_1},\ldots,v_{i_k}\}$ spans a simplex of $\tilde K$ if and only if
the complement $\{v_1,\ldots,v_m\}\setminus v_I$ does not span a simplex of
$K$. The product of cohomology classes of the complement of a coordinate
subspace arrangement was described in~\cite{dL} in combinatorial terms using
the complex $\tilde K$ and the above interpretation of the
Goresky--MacPherson formula.

In our paper we prefer to describe a coordinate subspace arrangement in terms
of the simplicial complex $K$ instead of $\tilde K$ because such an approach
reveals new connections between the topology of complements of
subspace arrangements, commutative algebra, and geometry of {\it toric
varieties}. We show that the complement $U(K)$ is homotopically equivalent to
what we call the {\it moment-angle complex} $\zk$ defined by the simplicial
complex $K$. This $\zk$ is a compact subset of a unit poly-disk
$(D^2)^m\subset\C^m$ invariant with respect to the standard $T^m$-action on
$(D^2)^m$.  At the same time $\zk$ is a homotopy fibre of cellular
embedding $i:\tbtk\hookrightarrow BT^m$, where $BT^m$ is the
$T^m$-classifying space with standard cellular structure, and $\tbtk$ is a
cell subcomplex whose cohomology is isomorphic to the {\it Stanley--Reisner
face ring\/} $\k(K)$ of simplicial complex $K$. Then we calculate the
cohomology algebra of $\zk$ (or $U(K)$) by means of the Eilenberg--Moore
spectral sequence. As the result, we obtain an algebraic description of the
cohomology algebra of $U(K)$ as the bigraded cohomology algebra
$\Tor_{\k[v_1,\ldots,v_m]}\bigl(\k(K),\k\bigr)$ of the face ring $\k(K)$. By
means of the standard Koszul resolution the latter can be expressed as the
cohomology of differential bigraded algebra
$\k(K)\otimes\Lambda[u_1,\ldots,u_m]$, where $\Lambda[u_1,\ldots,u_m]$ is an
exterior algebra, and the differential sends exterior generator $u_i$ to
$v_i\in\k(K)=\k[v_1,\ldots,v_m]/I$. The rational models of De Concini and
Procesi~\cite{dCP} and Yuzvinsky~\cite{Yu} also can be interpreted as an
application of the Koszul resolution to the cohomology of the complement a
subspace arrangement, however the role of the face ring became clear only
after our paper~\cite{BP2}.

If $K$ is an $(n-1)$-dimensional simplicial sphere (for instance, $K$ is the
boundary complex of an $n$-dimensional convex simplicial polytope), our
moment-angle complex $\zk$ turns to be a smooth $(m+n)$-dimensional manifold
(hence, $U(K)$ is homotopically equivalent to a smooth manifold). This
important particular case of our constructions was detailedly studied
in~\cite{BP1}, \cite{BP2}. Topological properties of the above manifolds
$\zk$ are of great interest because of their relations with combinatorics of
polytopes, symplectic geometry, and geometry of toric varieties; the
last thing was the starting point in our study of coordinate subspace
arrangements. The classical definition of toric varieties (see~\cite{Da},
\cite{Fu}) deals with the combinatorial object known as {\it fan}. However,
as it have been recently shown by several authors (see, for
example,~\cite{Au}, \cite{Ba}, \cite{Co}), in the case when the fan defining
a toric variety $M$ is simplicial, $M$ can be defined as the geometric
quotient of the complement $U(K)$ with respect to a certain action of the
algebraic torus $(\C^{*})^{m-n}$ (here $K$ is the simplicial complex defined
by the fan). Our moment-angle manifold $\zk$ is the pre-image of a regular
point in the image of the moment map $U(K)\to\R^{m-n}$ for the Hamiltonian
action of compact torus $T^{m-n}\subset(\C^{*})^{m-n}$.

In their paper~\cite{DJ} Davis and Januszkiewicz introduced the notion of
toric manifold (now also known as quasitoric manifold or unitary toric
manifold), which can be regarded as a natural topological extension of the
notion of smooth toric variety. A (quasi)toric manifold $M^{2n}$ admits a
smooth action of the torus $T^n$ that locally looks like the standard action
of $T^n$ on $\C^n$; the orbit space is required to be an $n$-dimensional
ball, invested with the combinatorial structure of a simple convex polytope
by the fixed point sets of appropriate subtori. Topology, geometry and
combinatorics of quasitoric manifolds are very beautiful; after the
pioneering paper~\cite{DJ} many new relations have been discovered by
different authors (see~\cite{BR1}, \cite{BR2}, \cite{BP1}, \cite{BP2},
\cite{Pa1}, \cite{Pa2}, and more references there). The dual complex to the
boundary complex of a simple polytope in the orbit space of a quasitoric
manifold is a simplicial sphere.  That is why many results from the present
paper may be considered as an extension of our previous constructions with
simplicial spheres to the case of general simplicial complex. We also mention
that some our definitions and constructions (such as the Borel construction
$B_TP$) firstly appeared in~\cite{DJ} in a different fashion; in this case we
have tried to retain initial notations.

The authors express special thanks to Nigel Ray for stimulating discussions
and fruitful collaboration which inspired some ideas and constructions from
this paper. We also grateful to Nataliya Dobrinskaya who have drawn our
attention to paper~\cite{Ba}, which reveals some connections between toric
varieties and coordinate subspace arrangements, and to Sergey Yuzvinsky who
informed us about the results of preprint~\cite{dL}.

\section{Homotopical realization of complement of a coordinate subspace
arrangement}

Let $\C^m$ be a complex $m$-dimensional space with coordinates
$z_1,\ldots,z_m$. For any index subset $I=\{i_1,\ldots,i_k\}$ denote by $L_I$
the $(m-k)$-dimensional coordinate subspace defined by the equations
$z_{i_1}=\ldots=z_{i_k}=0$. Note that $L_{\{1,\ldots,m\}}=\{0\}$ and
$L_{\varnothing}=\C^m$.

\begin{definition}
  A {\it coordinate subspace arrangement\/} $\mathcal A$ is a set
  of coordinate subspaces $L_I$. The {\it complement} of $\mathcal A$
  is the subset
  $$
    U(\mathcal A)=\C^m\setminus\bigcup_{L_I\in\mathcal A}L_I\subset\C^m.
  $$
\end{definition}
In the sequel we would distinguish the coordinate subspace arrangement
$\mathcal A$ regarded as an abstract set of subspaces and its support
$|\mathcal A|$~--- the subset
$\bigcup_{L_I\in\mathcal A}L_I\subset\C^m$. If $I\subset J$ and
$L_I\subset|\mathcal A|$, then $L_J\subset|\mathcal A|$.
If a coordinate subspace arrangement $\mathcal A$ contains a hyperplane
$z_i=0$, then its complement $U(\mathcal A)$ is represented as $U(\mathcal
A_0)\times\C^{*}$, where $\mathcal A_0$ is a coordinate subspace
arrangement in the hyperplane $\{z_i=0\}$, and $\C^{*}=\C\setminus\{0\}$.
Thus, for any coordinate subspace arrangement $\mathcal A$ the complement
$U(\mathcal A)$ decomposes as
$$
  U(\mathcal A)=U(\mathcal A')\times(\C^{*})^k,
$$
were $\mathcal A'$ is a coordinate arrangement in $\C^{m-k}$ that does not
contain any hyperplane. Keeping in mind this remark, we restrict ourself to
coordinate subspace arrangement without hyperplanes.

A coordinate
subspace arrangement $\mathcal A$ in $\C^m$ (without hyperplanes) defines a
simplicial complex $K(\mathcal A)$ with $m$ vertices $v_1,\ldots,v_m$ in the
following way: we say that a subset $v_I=\{v_{i_1},\ldots,v_{i_k}\}$ is a
$(k-1)$-simplex of $K(\mathcal A)$ if and only if $L_I\not\subset|\mathcal
A|$.

\begin{example}
  1) If $\mathcal A=\varnothing$, then $K(\mathcal A)$ is an
  $(m-1)$-dimensional simplex $\D^{m-1}$.

  2) If $\mathcal A=\{0\}$, then
  $K(\mathcal A)=\partial\D^{m-1}$ is the boundary of an $(m-1)$-simplex.
\end{example}
On the other hand, a simplicial complex $K$ on the vertex set
$\{v_1,\ldots,v_m\}$ defines an arrangement $\mathcal A(K)$ such that
$L_I\subset|\mathcal A|$ if and only if $v_I=\{v_{i_1},\ldots,v_{i_k}\}$ is
not a simplex of $K$. Note that if $K'\subset K$ is a subcomplex, then
$\mathcal A(K)\subset\mathcal A(K')$.  Thus, we have a reversing order
one-to-one correspondence between simplicial complexes on $m$ vertices and
coordinate subspace arrangements in $\C^m$ without hyperplanes.

Now let $U(K)=\C^m\setminus|\mathcal A(K)|$ denote the complement of the
coordinate subspace arrangement $\mathcal A(K)$.

\begin{example}
\label{uk}
1) If $K=\D^{m-1}$ is an $(m-1)$-simplex, then $U(K)=\C^m$.

2) If $K=\partial\D^{m-1}$, then $U(K)=\C^m\setminus\{0\}$.

3) If $K$ is a disjoint union of $m$ vertices, then $U(K)$ is obtained
by removing from $\C^m$ all codimension-two coordinate subspaces
$z_i=z_j=0$, $i,j=1,\ldots,m$.
\end{example}

Suppose that $\k$ is any field, which we refer to as the ground field.
Form a polynomial ring $\k[v_1,\ldots,v_m]$ where the $v_i$ are regarded as
indeterminates.

\begin{definition}
\label{fr}
  The {\it face ring\/} (or the {\it Stanley--Reisner ring\/}) $\k(K)$ of
  simplicial complex $K$ is the quotient ring $\k[v_1,\ldots,v_m]/I$, where
  $$
    I=\left(v_{i_1}\cdots v_{i_s}:
    \{v_{i_1},\ldots,v_{i_s}\}\; \text{ does not span a simplex in }K
    \right).
  $$
\end{definition}
Thus, the face ring is a quotient ring of polynomial ring by an ideal
generated by square free monomials of degree $\ge2$. We make $\k(K)$ a graded
ring by setting $\deg v_i=2$, $i=1,\ldots,m$.

\begin{example}
1) If $K=\D^{m-1}$, then $\k(K)=\k[v_1,\ldots,v_m]$.

2) If $K=\partial\D^{m-1}$ is the boundary complex of a
$(m-1)$-simplex, then $\k(K)=\k[v_1,\ldots,v_m]/(v_1\cdots v_m)$.
\end{example}

A compact torus $T^m$ acts
on $\C^m$ diagonally; since the arrangement $\mathcal A(K)$ consists of
coordinate subspaces, this action is also defined on $U(K)$. Denote by $B_TK$
the corresponding Borel construction:
\begin{equation}
\label{borel}
  B_TK=ET^m\times_{T^m}U(K),
\end{equation}
where $ET^m$ is the contractible space of universal $T^m$-bundle
$ET^m\to BT^m$ over the classifying space $BT^m=(\C P^{\infty})^m$. Thus,
$B_TK$ is the total space of bundle $B_TK\to BT^m$ with fibre $U(K)$.

The space $BT^m$ has a canonical cellular decomposition (that is, each
$\C P^\infty$ has one cell in each even dimension). For each index set
$I=\{i_1,\ldots,i_k\}$ one may consider the cellular subcomplex
$BT^k_I=BT^k_{i_1,\ldots,i_k}\subset BT^m$ homeomorphic to $BT^k$.

\begin{definition}
\label{tbtk}
  Given a simplicial complex $K$ with vertex set $\{v_1,\ldots,v_m\}$,
  define cellular subcomplex $\widetilde{B_TK}\subset BT^m$ as the union
  of $BT^k_I$ over all $I$ such that $v_I$ is a simplex of $K$.
\end{definition}

\begin{example}
\label{bouquet}
  Let $K$ be a disjoint union of $m$ vertices $v_1,\ldots,v_m$.
  Then $\widetilde{B_TK}$ is a bouquet of
  $m$ copies of $\C P^{\infty}$.
\end{example}

The cohomology ring of $BT^m$ is isomorphic to the polynomial ring
$\k[v_1,\ldots,v_m]$ (all cohomologies are with coefficients in the ground
field~$\k$).

\begin{lemma}
\label{cohomtbtk}
  The cohomology ring of $\widetilde{B_TK}$ is isomorphic to the face ring
  $\k(K)$. The embedding $i:\widetilde{B_TK}\hookrightarrow BT^m$
  induces the quotient epimorphism
  $i^{*}:\k[v_1,\ldots,v_m]\to\k(K)=\k[v_1,\ldots,v_m]/I$ in the cohomology.
\end{lemma}
\begin{proof}
The proof is by induction on the number of simplices of $K$.
If $K$ is a disjoint union of vertices $v_1,\ldots,v_m$, then
$\tbtk$ is a bouquet of $m$ copies of $\C P^{\infty}$ (see
Example~\ref{bouquet}). In degree zero $H^{*}(\tbtk)$ is just $\k$, while in
degrees $\ge 1$ it is isomorphic to $\k[v_1]\oplus\cdots\oplus\k[v_m]$.
Therefore, $H^{*}(\tbtk)=k[v_1,\ldots,v_m]/I$, where $I$ is the ideal
generated by all square free monomials of degree $\ge 2$, and $i^{*}$ is the
projection onto the quotient ring. Thus, the lemma holds for $\dim K=0$.

Now suppose that the simplicial complex $K$ is obtained from the simplicial
complex $K'$ by adding one $(k-1)$-dimensional simplex
$v_I=\{v_{i_1},\ldots,v_{i_k}\}$. By the inductive hypothesis, the
lemma holds for $K'$, that is,
$i^{*}H^{*}(BT^m)=H^{*}(\widetilde{B_TK'})=\k(K')=\k[v_1,\ldots,v_m]/I'$.
By Definition~\ref{tbtk}, $\tbtk$ is obtained from $\widetilde{B_TK'}$ by
adding the subcomplex $BT^k_{i_1,\ldots,i_k}\subset BT^m$.
Then $H^{*}(\widetilde{B_TK'}\cup
BT^k_{i_1,\ldots,i_k})=\k[v_1,\ldots,v_m]/I=\k(K'\cup v_I)$, where $I$
is generated by $I'$ and $v_{i_1}v_{i_2}\cdots v_{i_k}$.
\end{proof}

Let $I^m$ be the standard $m$-dimensional cube in $\R^m$:
$$
  I^m=\{(y_1,\ldots,y_m)\in\R^m\,:\;0\le y_i\le 1,\,i=1,\ldots,m\}.
$$
A simplicial complex $K$ with $m$ vertices $v_1,\ldots,v_m$ defines a cubical
complex $\mathcal C_K$ embedded canonically into the boundary complex of
$I^m$ in the following way:

\begin{definition}
\label{cubcom}
  For each $(k-1)$-dimensional simplex $v_J=\{v_{j_1},\ldots,v_{j_k}\}$ of
  $K$ denote by $C_J$ the $k$-dimensional face of $I^m$ defined
  by $m-k$ equations
  $$
    y_i=1,\quad i\notin \{j_1,\ldots,j_k\}.
  $$
  Then define cubical subcomplex $\mathcal C_K\subset I^m$ as the union of
  $C_J$ over all simplices $v_J$ of $K$.
\end{definition}

\begin{remark}
  Our cubical subcomplex $\mathcal C_K\subset I^m$ is a geometrical
  realization of an abstract cubical complex in the cone over the
  barycentric subdivision of $K$ (see~\cite[p.~434]{DJ}).
  Indeed, let $\D^{m-1}$ be an $(m-1)$-dimensional simplex on the vertex
  set $\{v_1,\ldots,v_m\}$, and $\hat\D^{m-1}$ a barycentric subdivision of
  $\D^{m-1}$, that is, $\hat\D^{m-1}$ has a vertex for each simplex $v_J$ of
  $\D^{m-1}$. Construct a map $\iota:\hat\D^{m-1}\to I^m$ by sending vertex
  $v_J$ of $\hat\D^{m-1}$ to the vertex of $I^m$ having coordinates $y_j=0$
  for $j\in J$ and $y_j=1$ for $j\notin J$ and then extending this map
  linearly on each simplex of $\hat\D^{m-1}$. The image of $\hat\D^{m-1}$
  under the constructed map is the union of faces of $I^m$ meeting at zero.
  Then build a map $C\iota$ from the cone $C\hat\D^{m-1}$ over $\hat\D^{m-1}$
  to $I^m$ by sending the vertex of the cone to the vertex $(1,\ldots,1)$ of
  the cube and extending linearly on simplices of $C\hat\D^{m-1}$. The image
  of $C\hat\D^{m-1}$ under the map $C\iota$ is the whole cube $I^m$. Now let
  $K$ be a simplicial complex on the vertex set $\{v_1,\ldots,v_m\}$. Once a
  numeration of vertices is fixed, we may view $K$ as a simplicial subcomplex
  of $\D^{m-1}$. Then our cubical complex $\mathcal C_K\subset I^m$ from
  Definition~\ref{cubcom} is nothing but the image $C\iota(C\hat K)$ of the
  cone over the barycentric subdivision of $K$ under the map $C\iota$.
\end{remark}

\begin{example}
  The cubical complex $\mathcal C_K$ in the cases when $K$ is a disjoint
  union of 3 vertices and the boundary complex of a 2-simplex is indicated on
  Figure~1~a) and~b) respectively.
  \begin{figure}
  \begin{picture}(120,45)
  \put(15,5){\line(1,0){25}}
  \put(15,5){\line(0,1){25}}
  \put(40,5){\line(1,1){10}}
  \put(50,15){\line(0,1){25}}
  \put(50,40){\line(-1,0){25}}
  \put(25,40){\line(-1,-1){10}}
  \put(40,5){\circle*{2}}
  \put(40,30){\circle*{2}}
  \put(15,30){\circle*{2}}
  \put(50,40){\circle*{2}}
  \multiput(40,29.3)(0,0.1){16}{\line(1,1){10}}
  \multiput(25,15)(5,0){5}{\line(1,0){3}}
  \multiput(25,15)(0,5){5}{\line(0,1){3}}
  \multiput(25,15)(-5,-5){2}{\line(-1,-1){4}}
  \put(26,16){$0$}
  \put(30,0){a)}
  \put(75,5){\line(1,0){25}}
  \put(75,5){\line(0,1){25}}
  \put(100,5){\line(1,1){10}}
  \put(110,15){\line(0,1){25}}
  \put(110,40){\line(-1,0){25}}
  \put(85,40){\line(-1,-1){10}}
  \put(100,5){\circle*{2}}
  \put(110,15){\circle*{2}}
  \put(100,30){\circle*{2}}
  \put(75,30){\circle*{2}}
  \put(85,40){\circle*{2}}
  \put(75,5){\circle*{2}}
  \put(110,40){\circle*{2}}
  \multiput(100,29.3)(0,0.1){16}{\line(1,1){10}}
  \multiput(75,29.3)(0,0.1){16}{\line(1,1){10}}
  \multiput(100,4.3)(0,0.1){16}{\line(1,1){10}}
  \multiput(85,15)(5,0){5}{\line(1,0){3}}
  \multiput(85,15)(0,5){5}{\line(0,1){3}}
  \multiput(85,15)(-5,-5){2}{\line(-1,-1){4}}
  \multiput(78,5)(4,0){6}{\line(0,1){25}}
  \multiput(78,30)(4,0){6}{\line(1,1){10}}
  \multiput(100,8)(0,4){6}{\line(1,1){10}}
  \put(86,16){$0$}
  \put(90,0){b)}
  \linethickness{1mm}
  \put(15,30){\line(1,0){25}}
  \put(40,5){\line(0,1){25}}
  \put(75,30){\line(1,0){25}}
  \put(75,5){\line(1,0){25}}
  \put(85,40){\line(1,0){25}}
  \put(100,5){\line(0,1){25}}
  \put(75,5){\line(0,1){25}}
  \put(110,15){\line(0,1){25}}
  \end{picture}
  \caption{The cubical complex $\mathcal C_K$.}
  \end{figure}
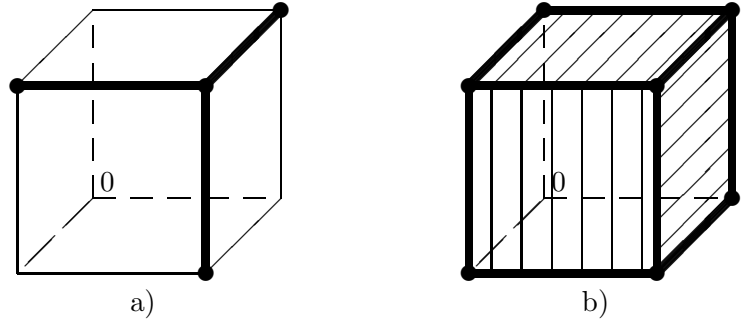
\end{example}

\begin{remark}
  In the case when $K$ is the dual to the boundary complex of an
  $n$-dimensional simple polytope $P^n$, the cubical complex $\mathcal C_K$
  coincides with the cubical subdivision of $P^n$ studied in~\cite{BP2}.
\end{remark}

The orbit space of the diagonal action of $T^m$ on $\C^m$ is the positive
cone
$$
  \R^m_+=\{(y_1,\ldots,y_m)\in\R^m\,:\;y_i\ge 0,\,i=1,\ldots,m\}.
$$
The orbit map $\C^m\to\R^m_+$ can be given by
$(z_1,\ldots,z_m)\to(|z_1|^2,\ldots,|z_m|^2)$. If we restrict the above
action to the standard poly-disk
$$
 (D^2)^m=\{(z_1,\ldots,z_m)\in\C^m:\;\:|z_i|\le 1,\,i=1,\ldots,m\}\subset\C^m,
$$
then the corresponding orbit space would be the standard cube
$I^m\subset\R^m_+$.

Let $U_{\R}(K)\subset\R^m_+$ denote the orbit space $U(K)/T^m$. Note that if
we regard $\R^m_+$ as a subset in $\C^m$, then $U_{\R}(K)$ is the ``real
part": $U_{\R}(K)=U(K)\cap\R^m_+$.

\begin{definition}
\label{zk}
  The {\it equivariant moment-angle complex} $\zk\subset\C^m$ corresponding
  to a simplicial complex $K$ is the $T^m$-space defined
  from the commutative diagram
  $$
    \begin{CD}
      \zk @>>> (D^2)^m\\
      @VVV @VVV\\
      \mathcal C_K @>>> I^m,
    \end{CD}
  $$
  where the right vertical arrow denotes the orbit map for the diagonal action
  of $T^m$, and the lower horizontal arrow denotes the embedding of the cubical
  complex $\mathcal C_K$ to $I^m$.
\end{definition}

\begin{lemma}
\label{zu}
  $\mathcal C_K\subset U_{\R}(K)$ and $\zk\subset U(K)$.
\end{lemma}
\begin{proof}
Definition~\ref{zk} shows that the second assertion follows from the first
one. To prove the first assertion we mention that if a point
$a=(y_1,\ldots,y_m)\in\mathcal C_K$ has $y_{i_1}=\ldots=y_{i_k}=0$, then
$v_I=\{v_{i_1},\ldots,v_{i_k}\}$ is a simplex of $K$, hence
$L_I\not\subset\mathcal A(K)$.
\end{proof}

\begin{lemma}
\label{he1}
  $U(K)$ is equivariantly
  homotopy equivalent to the moment-angle complex~$\zk$.
\end{lemma}
\begin{proof}
We construct a retraction $r:U_{\R}(K)\to\mathcal C_K$ that is covered by
an equivariant retraction $U(K)\to\zk$. The latter would be a required
homotopy equivalence.

The retraction $r:U_{\R}(K)\to\mathcal C_K$ is constructed inductively. We
start from the boundary complex of an $(m-1)$-simplex and remove simplices of
positive dimensions until we obtain $K$. On each step we construct a
retraction, and the composite map would be required retraction $r$. If
$K=\partial\D^{m-1}$ is the boundary complex of an $(m-1)$-simplex, then
$U_{\R}(K)=\R^m_+\setminus\{0\}$ and the retraction $r$ is shown on Figure~2.
\begin{figure}
  \begin{picture}(120,45)
  \put(45,5){\circle{2}}
  \put(80,5){\circle*{2}}
  \put(45,40){\circle*{2}}
  \put(80,40){\circle*{2}}
  \put(46,5){\vector(1,0){24}}
  \put(46,5){\line(1,0){34}}
  \put(45.8,5.2){\vector(4,1){24}}
  \put(45.8,5.2){\line(4,1){34}}
  \put(45.5,5.5){\vector(2,1){24}}
  \put(45.5,5.5){\line(2,1){34}}
  \put(45.8,5.8){\vector(4,3){24}}
  \put(45.8,5.8){\line(4,3){34}}
  \put(46,6){\vector(1,1){24}}
  \put(46,6){\line(1,1){34}}
  \put(45,6){\vector(0,1){24}}
  \put(45,6){\line(0,1){34}}
  \put(45.2,5.8){\vector(1,4){6}}
  \put(45.2,5.8){\line(1,4){8.5}}
  \put(45.5,5.5){\vector(1,2){12}}
  \put(45.5,5.5){\line(1,2){17}}
  \put(45.8,5.2){\vector(3,4){18}}
  \put(45.8,5.2){\line(3,4){26}}
  \linethickness{1mm}
  \put(80,5){\line(0,1){35}}
  \put(45,40){\line(1,0){35}}
    \end{picture}
  \caption{The retraction $r:U_{\R}(K)\to\mathcal C_K$ for
  $K=\partial\D^{m-1}$.}
\end{figure}
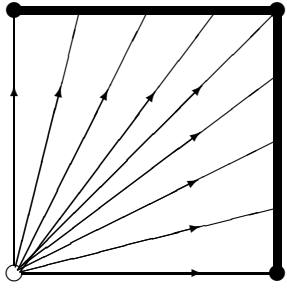
Now suppose that the simplicial complex $K$ is obtained by removing one
$(k-1)$-dimensional simplex $v_J=\{v_{j_1},\ldots,v_{j_k}\}$ from simplicial
complex $K'$. By the inductive hypothesis, the
lemma holds for $K'$, that is, there is a retraction
$r':U_{\R}(K')\to\mathcal C_{K'}$ with the required properties. Let us consider
the face $C_J\subset I^m$ (see Definition~\ref{cubcom}). Since
$v_J$ is not a simplex of $K$, the point $a$ having coordinates
$y_{j_1}=\ldots=y_{j_k}=0$, $y_i=1$, $i\notin\{j_1,\ldots,j_k\}$, do not
belong to $U(K)$. Hence, we may apply the retraction from Figure~2 on the
face $C_J$, starting from the point~$a$. Denote this retraction by $r_J$. Now
take $r=r_J\circ r'$. It is easy to see that this $r$ is exactly the required
retraction.
\end{proof}

\begin{example}
\label{zp}
1) If $K=\partial\D^{m-1}$ is the boundary complex of an $(m-1)$-simplex,
  then $\zk$ is homeomorphic to $(2m-1)$-dimensional sphere $S^{2m-1}$.

2) If $K$ is the dual to the boundary complex of a
  $n$-dimensional simple polytope $P^n$, then $\zk$ is homeomorphic to a
  smooth $(m+n)$-dimensional manifold. This manifold, denoted $\mathcal Z_P$,
  is the main object of study in~\cite{BP2}.
\end{example}

\begin{corollary}
  The Borel construction $ET^m\times_{T^m}\zk$ is homotopy equivalent to
  $B_TK$.
\end{corollary}
\begin{proof}
The retraction $r:U(K)\to\zk$ constructed in the proof
of Lemma~\ref{he1} is equivariant with respect to the $T^m$-actions on $U(K)$
and $\zk$. Since $B_TK=ET^m\times_{T^m}U(K)$, the corollary follows.
\end{proof}

In what follows we do not distinguish the Borel constructions
$ET^m\times_{T^m}\zk$ and $B_TK=ET^m\times_{T^m}U(K)$.

\begin{theorem}
\label{he2}
  The cellular embedding $i:\tbtk\hookrightarrow BT^m$ (see
  {\rm Definition~\ref{tbtk}}) and the fibration $p:B_TK\to BT^m$
  (see~{\rm (\ref{borel})}) are homotopically equivalent. In particular,
  $\tbtk$ and $B_TK$ are of same homotopy type.
\end{theorem}
\begin{proof}
Let $\pi:\zk\to\mathcal C_K$ denote the orbit map for the torus action on
the moment-angle complex $\zk$ (see Definition~\ref{zk}). For each subset
$I=\{i_1,\ldots,i_k\}\subset\{1,\ldots,m\}$ denote by $B_I$ the following
subset of the poly-disk $(D^2)^m$: $B_I=B_1\times\cdots\times B_m\subset
D^2\times\cdots\times D^2=(D^2)^m$, where $B_i$ is the disk $D^2$ if $i\in
I$, and $B_i$ is the boundary $S^1$ of $D^2$ if $i\notin I$. Thus, $B_I\cong
(D^2)^k\times T^{m-k}$, where $k=|I|$. It is easy to see that if $C_I$ is a
face of cubical complex $\mathcal C_K$ (see Definition~\ref{cubcom}) then
$\pi^{-1}(C_I)=B_I$. Since for $I\subset J$ the $B_I$ is canonically
identified with a subset of $B_J$, we see that those $B_I$ for which $v_I$ is
a simplex of $K$ fit together to yield $\zk$. (This idea can be used to prove
that $\zk$ is a smooth manifold provided that $K$ is the dual to the boundary
complex of a simple polytope, see~\cite[Theorem~2.4]{BP2}.)

For any simplex $v_I\subset K$ the subset $B_I\subset\zk$ is invariant with
respect to the $T^m$-action on $\zk$. Hence, the Borel construction
$B_TK=ET^m\times_{T^m}\zk$ is patched from Borel
constructions $ET^m\times_{T^m}B_I$ (compare this with the local construction
of $B_TP$ from~\cite[p.~435]{DJ}). The latter can be factorized as
$ET^m\times_{T^m}B_I=\bigl(ET^k\times_{T^k}(D^2)^k\bigr)\times ET^{m-n}$,
which is homotopically equivalent to $BT^k_I$. Hence, the restriction of the
projection $p:B_TK\to BT^m$ to $ET^m\times_{T^m}B_I$ is homotopically
equivalent to the embedding $BT^k_I\hookrightarrow BT^m$. These homotopy
equivalences for all simplices $v_I\subset K$ fit together to yield a
required homotopy equivalence between $p:B_TK\to BT^m$ and
$i:\tbtk\hookrightarrow BT^m$.
\end{proof}

\begin{corollary}
\label{homfib}
  The complement $U(K)$ of a coordinate subspace arrangement is a homotopy
  fibre of the cellular embedding $i:\tbtk\hookrightarrow BT^m$.\ep
\end{corollary}

\begin{corollary}
\label{eqv}
  The $T^m$-equivariant cohomology ring $H^{*}_{T^m}\bigl(U(K)\bigr)$ is
  isomorphic to the face ring $\k(K)$.
\end{corollary}
\begin{proof}
  We have $H^{*}_{T^m}\bigl(U(K)\bigr)=H^{*}\bigl(ET^m\times_{T^m}U(K)\bigr)=
  H^{*}(B_TK)$. Now, the corollary follows from
  Lemma~\ref{cohomtbtk} and Theorem~\ref{he2}.
\end{proof}

\section{Cohomology ring of $U(K)$}

Suppose that we are given a $\k[v_1,\ldots,v_m]$-free resolution of the face
ring $\k(K)$ as a graded
module over the polynomial ring $\k[v_1,\ldots,v_m]$:
\begin{equation}
\label{resol}
  \begin{CD}
  0\to R^{-h} @>d^{-h}>> R^{-h+1} @>d^{-h+1}>> \cdots\to R^{-1}
  \stackrel{d^{-1}}{\longrightarrow} R^0
  \stackrel{d^{0}}{\longrightarrow} \k(K)\to 0
  \end{CD}
\end{equation}
(note that the Hilbert syzygy theorem shows that $h\le m$ above). Applying
the functor $\otimes_{\k[v_1,\ldots,v_m]}\k$ to~(\ref{resol}) we obtain a
cochain complex:
$$
  \begin{CD}
   0 \longrightarrow R^{-h}\otimes_{\k[v_1,\ldots,v_m]}\k \longrightarrow
   \cdots \longrightarrow R^0\otimes_{\k[v_1,\ldots,v_m]}\k
   \longrightarrow 0,
  \end{CD}
$$
whose cohomology modules are denoted
$\Tor^{-i}_{\k[v_1,\ldots,v_m]}\bigl(\k(K),\k\bigr)$. Since all $R^{-i}$
in~(\ref{resol}) are graded $\k[v_1,\ldots,v_m]$-modules,
$\Tor^{-i}_{\k[v_1,\ldots,v_m]}\bigl(\k(K),\k\bigr)=\bigoplus_j
\Tor^{-i,j}_{\k[v_1,\ldots,v_m]}\bigl(\k(K),\k\bigr)$ is a graded $\k$-module,
and
\begin{equation}
\label{tor}
  \Tor_{\k[v_1,\ldots,v_m]}\bigl(\k(K),\k\bigr)=\bigoplus_{i,j}
  \Tor^{-i,j}_{\k[v_1,\ldots,v_m]}\bigl(\k(K),\k\bigr)
\end{equation}
becomes a {\it bigraded} $\k$-module. Note that its non-zero elements have
non-positive first grading and non-negative even
second grading (since $\deg v_i$=2).
The bigraded $\k$-module~(\ref{tor}) can be also regarded as a
one-graded module with respect to the total degree $-i+j$. The Betti numbers
\begin{align*}
  \beta^{-i}\bigl(\k(K)\bigr)&=
  \dim_{\k}\Tor^{-i}_{\k[v_1,\ldots,v_m]}\bigl(\k(K),\k\bigr)\\
  \intertext{and}
  \beta^{-i,2j}\bigl(\k(K)\bigr)&=
  \dim_{\k}\Tor^{-i,2j}_{\k[v_1,\ldots,v_m]}\bigl(\k(K),\k\bigr)
\end{align*}
are of great
interest in geometric combinatorics; they were studied by different authors
(see, for example,~\cite{St}). We mention only one theorem due to Hochster,
which reduces calculation of $\beta^{-i,2j}\bigl(\k(K)\bigr)$ to calculating
the homology of subcomplexes of~$K$.
\begin{theorem}[\rm Hochster~\cite{Ho}, \cite{St}]
\label{hoch}
  The Hilbert series $\sum_{j}\beta^{-i,2j}\bigl(\k(K)\bigr)t^{2j}$
  of $\Tor^{-i}_{\k[v_1,\ldots,v_m]}\bigl(\k(K),\k\bigr) $
  can be calculated as
  $$
    \sum_{j}\beta^{-i,2j}\bigl(\k(K)\bigr)t^{2j}=
    \sum_{I\subset\{v_1,\ldots,v_m\}}
    \bigl(\dim_{\k}\tilde{H}_{|I|-i-1}(K_I)\bigr)t^{2|I|},
  $$
  where $K_I$ is the subcomplex of $K$ consisting of all simplices with
  vertices in~$I$.\ep
\end{theorem}
Note that calculation of $\beta^{-i,2j}\bigl(\k(K)\bigr)$ using this theorem
is very involved even for small~$K$.

It turns out that $\Tor_{\k[v_1,\ldots,v_m]}\bigl(\k(K),\k\bigr)$ is a
bigraded {\it algebra} in a natural way, and the associated one-graded
algebra is exactly $H^{*}\bigl(U(K)\bigr)$:
\begin{theorem}
\label{cohom1}
  The following isomorphism of graded algebras holds:
  $$
   H^{*}\bigl(U(K)\bigr)\cong\Tor_{\k[v_1,\ldots,v_m]}\bigl(\k(K),\k\bigr)
  $$
\end{theorem}
\begin{proof}
Let us consider the commutative diagram
\begin{equation}
\label{wuk}
  \begin{CD}
    \widetilde{U(K)} @>>> ET^m\\
    @VVV @VVV\\
    \tbtk @>i>> BT^m,
  \end{CD}
\end{equation}
where the left vertical arrow is the induced fibre bundle.
Corollary~\ref{homfib} shows that $\widetilde{U(K)}$ is homotopically
equivalent to $U(K)$.

From~(\ref{wuk}) we obtain that the cellular
cochain algebras $C^{*}(\tbtk)$ and
$C^{*}(ET^m)$ are modules over $C^{*}(BT^m)$. It is clear from the proof of
Lemma~\ref{cohomtbtk} that $C^{*}(\tbtk)=\k(K)$ and
$i^{*}:C^{*}(BT^m)=\k[v_1,\ldots,v_m]\to\k(K)=C^{*}(\tbtk)$ is the quotient
epimorphism. Since $ET^m$ is contractible, we have a chain equivalence
$C^{*}(ET^m)\to\k$. Therefore, there is an isomorphism
\begin{equation}
\label{chain}
 \Tor_{C^{*}(BT^m)}\bigl(C^{*}(\tbtk),C^{*}(ET^m)\bigr)\cong
 \Tor_{\k[v_1,\ldots,v_m]}\bigl(\k(K),\k\bigr).
\end{equation}

The Eilenberg--Moore spectral sequence (see~\cite[Theorem~1.2]{Sm})
of commutative square~(\ref{wuk}) has
the $E_2$-term
$$
  E_2=\Tor_{H^{*}(BT^m)}\bigl(H^{*}(\tbtk),H^{*}(ET^m)\bigr)
$$
and converges to $\Tor_{C^{*}(BT^m)}(C^{*}(\tbtk),C^{*}(ET^m))$. Since
$$
  \Tor_{H^{*}(BT^m)}(H^{*}(\tbtk),H^{*}(ET^m))=
  \Tor_{\k[v_1,\ldots,v_m]}\bigl(\k(K),\k\bigr),
$$
it follows from~(\ref{chain})
that the spectral sequence collapses at the $E_2$ term, that is,
$E_2=E_{\infty}$. Now, Proposition~3.2 of~\cite{Sm} shows that the module
$\Tor_{C^{*}(BT^m)}\bigl(C^{*}(\tbtk),C^{*}(ET^m)\bigr)$ is an algebra
isomorphic to $H^{*}\bigl(\widetilde{U(K)}\bigr)$, which concludes the proof.
\end{proof}

Our next theorem gives an explicit description of the algebra
$H^{*}\bigl(U(K)\bigr)$ as the cohomology algebra of
a simple differential bigraded algebra. We consider the tensor product
$\k(K)\otimes\Lambda[u_1,\ldots,u_m]$ of the face ring
$\k(K)=\k[v_1,\ldots,v_m]/I$ and an exterior algebra
$\Lambda[u_1,\ldots,u_m]$
on $m$ generators and make it a differential bigraded algebra by setting
\begin{gather}
  \bideg v_i=(0,2),\quad\bideg u_i=(-1,2),\notag\\
  \label{diff}
  d(1\otimes u_i)=v_i\otimes 1,\quad d(v_i\otimes 1)=0,
\end{gather}
and requiring that $d$ be a derivation of algebras.

\begin{theorem}
\label{cohom2}
  The following isomorphism of graded algebras holds:
  $$
    H^{*}\bigl(U(K)\bigr)\cong
    H\bigl[\k(K)\otimes\Lambda[u_1,\ldots,u_m],d\bigr],
  $$
  where in the right hand side stands the one-graded algebra associated to
  the bigraded cohomology algebra.
\end{theorem}
\begin{proof}
One can make $\k$ a $\k[v_1,\ldots,v_m]$-module by means of the homomorphism
that sends 1 to 1 and $v_i$ to 0.
Let us consider the Koszul resolution
(see, for example~\cite[Chapter~VII, \S~2]{Ma}) of $\k$
regarded as a $\k[v_1,\ldots,v_m]$-module:
$$
  \bigl[\k[v_1,\ldots,v_m]\otimes\Lambda[u_1,\ldots,u_m],d\bigr],
$$
where the differential $d$
is defined as in~(\ref{diff}). Since the bigraded torsion product
$\Tor_{\k[v_1,\ldots,v_m]}(\ ,\ )$ is a symmetric function of its arguments,
one has
$$
  \Tor_{\Gamma}\bigl(\k(K),\k\bigr)\\
  =H\bigl[\k(K)\otimes_{\Gamma}\Gamma\otimes
  \Lambda[u_1,\ldots,u_n],d\bigr]\\
  =\bigl[\Gamma\otimes\Lambda[u_1,\ldots,u_m],d\bigr],
$$
where we denoted $\Gamma=\k[v_1,\ldots,v_m]$. Since
$H^{*}\bigl(U(K)\bigr)\cong\Tor_{\Gamma}\bigl(\k(K),\k\bigr)$ by
Theorem~\ref{cohom1}, we obtain the required isomorphism
\end{proof}

Note that the above theorem not only calculates the cohomology algebra of
$U(K)$, but also makes this algebra {\it bigraded}.

\begin{corollary}
\label{ls}
  The Leray--Serre spectral sequence of the bundle
  $\widetilde{U(K)}\to\tbtk$ with fibre $T^m$
  (see~{\rm(\ref{wuk})}) collapses at the $E_3$ term.
\end{corollary}
\begin{proof}
The spectral sequence under consideration converges to
$H^{*}\bigl(\widetilde{U(K)}\bigr)=H^{*}\bigl(U(K)\bigr)$ and has
$$
  E_2=H^{*}(\tbtk)\otimes H^{*}(T^m)=
  \k(K)\otimes\Lambda[u_1,\ldots,u_m].
$$
It is easy to see that the differential in the $E_2$ term acts as
in~(\ref{diff}). Hence, $E_3=H[E_2,d]=
H\bigl[\k(K)\otimes\Lambda[u_1,\ldots,u_m]\bigr]=H^{*}\bigl(U(K)\bigr)$ by
Theorem~\ref{cohom2}.
\end{proof}

\begin{proposition}
\label{monom}
  Suppose that a monomial
  $$
  v^{\alpha_1}_{i_1}\ldots v^{\alpha_p}_{i_p} u_{j_1}\ldots
  u_{j_q}\in\k(K)\otimes\Lambda[u_1,\ldots,u_m],
  $$
  where $i_1<\ldots<i_p$, $j_1<\ldots<i_q$, represents a non-trivial
  cohomology class
  in $H^{*}\bigl(U(K)\bigr)$. Then $\alpha_1=\ldots=\alpha_p=1$,
  $\{v_{i_1},\ldots,v_{i_p}\}$ spans a simplex of $K$, and
  $\{i_1,\ldots,i_p\}\cap\{j_1,\ldots,j_q\}=\varnothing$.
\end{proposition}
\begin{proof}
See~\cite[Lemma~5.3]{BP2}.
\end{proof}

As it was mentioned above (see Example~\ref{zp}), if $K$ is the boundary
complex of a convex simplicial polytope (or, equivalently, $K$ is the dual to
the boundary complex of a simple polytope) or at least a simplicial sphere,
then $U(K)$ has homotopy type of a smooth manifold $\zk$. It was shown
in~\cite[Theorem~2.10]{BP2} that the corresponding homotopy equivalence can
be interpreted as the orbit map $U(K)\to U(K)/\R^{m-n}\cong\zk$ with respect
to a certain action of $\R^{m-n}$ on $U(K)$.

The coordinate subspace arrangement $\mathcal A(K)$ and its complement $U(K)$
play important role in the theory of toric varieties and symplectic geometry
(see, for example,~\cite{Au}, \cite{Ba}, \cite{Co}). More precisely, any
$n$-dimensional simplicial toric variety $M$ defined by a (simplicial) fan
$\Sigma$ in $\Z^n$ with $m$ one-dimensional cones can be obtained as the
geometric quotient $U(K_{\Sigma})/G$. Here $G$ is a subgroup of the complex
torus $(\C^{*})^{m}$ isomorphic to $(\C^{*})^{m-n}$ and $K_{\Sigma}$ is
the simplicial complex defined by the fan $\Sigma$
($i$-simplices of $K_{\Sigma}$ correspond to $(i+1)$-dimensional cones of
$\Sigma$). A smooth projective toric variety $M$
is a symplectic manifold of real dimension $2n$.
This manifold can be constructed
by the process of {\it symplectic reduction} in the following way. Let
$G_{\R}\cong T^{m-n}$ denote the maximal compact subgroup of $G$, and let
$\mu:\C^m\to\R^{m-n}$ be the moment map for the Hamiltonian action of
$G_{\R}$ on $\C^m$. Then for each regular value $a\in\R^{m-n}$ of
$\mu$ there is a diffeomorphism
$$
  \mu^{-1}(a)/G_{\R}\longrightarrow U(K_{\Sigma})/G=M
$$
(see~\cite{Co} for more information). In this situation
it can be easily seen that $\mu^{-1}(a)$ is exactly our manifold $\zk$ for
$K=K_{\Sigma}$.

\begin{example}
  Let $G\cong\C^{*}$ be the diagonal subgroup in $(\C^{*})^{n+1}$ and
  $K_{\Sigma}$ be the boundary complex of an $n$-simplex. Then
  $U(K_{\Sigma})=\C^{n+1}\setminus\{0\}$ and
  $M=\C^{n+1}\setminus\{0\}/\C^{*}$ is the complex projective space $\C P^n$.
  The moment map $\mu:\C^m\to\R$ takes $(z_1,\ldots,z_m)\in\C^m$ to
  $\frac12(|z_1|^2+\ldots+|z_m|^2)$ and for $a\ne0$ one has $\mu^{-1}(a)\cong
  S^{2n+1}\cong\zk$ (see Example~\ref{zp}).
\end{example}

In the case when $K$ is a simplicial sphere (hence, the complement $U(K)$ is
homotopically equivalent to the smooth manifold $\zk$), there is Poincar\'e
duality defined in the cohomology ring of $U(K)$.

\begin{proposition}
\label{poin}
  Suppose that $K$ is a simplicial sphere of dimension $n-1$, hence,
  $U(K)$ is homotopically equivalent to the smooth manifold $\zk$. Then

  1) The Poincar\'e duality in $H^{*}\bigl(U(K)\bigr)$ regards the bigraded
  structure defined by {\rm Theorem~\ref{cohom2}}. More precisely, if
  $\alpha\in H^{-i,2j}(U(K)\bigr)$ is a cohomology class, then its Poincar\'e
  dual $D\alpha$ belongs to $H^{-(m-n)+i,2(m-j)}$.

  2) Let $\{v_{i_1},\ldots,v_{i_n}\}$ be an $(n-1)$-simplex of $K$ and let
  $j_1<\ldots<j_{m-n}$,
  $\{i_1,\ldots,i_n,j_1,\ldots,j_{m-n}\}=\{1,\ldots,m\}$.
  Then the value of the element
  $$
    v_{i_1}\cdots v_{i_n}u_{j_1}\cdots
    u_{j_{m-n}}\in H^{m+n}\bigl(U(K)\bigr)\cong H^{m+n}(\zk)
  $$
  on the fundamental class of $\zk$ equals $\pm1$.

  3) Let $\{v_{i_1},\ldots,v_{i_n}\}$ and
  $\{v_{i_1},\ldots,v_{i_{n-1}},v_{j_1}\}$ be two $(n-1)$-simplices of $K$
  having common $(n-2)$-face $\{v_{i_1},\ldots,v_{i_{n-1}}\}$, and
  $j_1,\ldots,j_{m-n}$ be as in~2). Then
  $$
    v_{i_1}\cdots v_{i_n}u_{j_1}\cdots u_{j_{m-n}}=
    v_{i_1}\cdots v_{i_{n-1}}v_{j_1}u_{i_n}u_{j_2}\cdots u_{j_{m-n}}
  $$
  in $H^{m+n}\bigl(U(K)\bigr)$.
\end{proposition}
\begin{proof}
For the proof of 1) and 2) see~\cite[Lemma~5.1]{BP2}. To prove 3) we just
mention that
\begin{multline*}
  d(v_{i_1}\cdots v_{i_{n-1}}u_{i_n}u_{j_1}u_{j_2}\cdots u_{j_{m-n}})\\
  =v_{i_1}\cdots v_{i_n}u_{j_1}\cdots u_{j_{m-n}}-
  v_{i_1}\cdots v_{i_{n-1}}v_{j_1}u_{i_n}u_{j_2}\cdots u_{j_{m-n}}
\end{multline*}
in $\k(K)\otimes\Lambda[u_1,\ldots,u_m]$ (see~(\ref{diff})).
\end{proof}

A simplicial complex $K$ is called {\it Cohen--Macaulay}, if its face ring
$\k(K)$ is a Cohen--Macaulay algebra, that is, $\k(K)$ is a
finite-dimensional free module over a polynomial ring $\k[t_1,\ldots,t_n]$
(here $n$ is the maximal number of algebraically independent elements of
$\k(K)$). Equivalently, $\k(K)$ is a Cohen--Macaulay algebra if it admits a
{\it regular sequence} $\{\l_1,\ldots,\l_n\}$, that is, a set of $n$
homogeneous elements such that $\l_{i+1}$ is not a zero divisor in
$\k(K)/(\l_1,\ldots,\l_i)$ for $i=0,\ldots,n-1$. If $K$ is a Cohen--Macaulay
complex and $\k$ is of infinite characteristic, then $\k(K)$ admits a regular
sequence of degree-two elements (remember that we set $\deg v_i=2$ in
$\k(K)$), that is, $\l_i=\l_{i1}v_1+\l_{i2}v_2+\ldots+\l_{im}v_m$,
$i=1,\ldots,n$.

\begin{theorem}
\label{cohom3}
  Suppose that $K$ is a Cohen--Macaulay complex and $J=(\l_1,\ldots,\l_n)$ is
  an ideal in $\k(K)$ generated by a regular sequence. Then the following
  isomorphism of bigraded algebras holds:
  $$
    H^{*}\bigl(U(K)\bigr)\cong
    H\bigl[\k(K)/\!J\otimes\Lambda[u_1,\ldots,u_{m-n}],d\bigr],
  $$
  where the gradings and differential in the right hand side are defined as
  follows:
  \begin{align*}
    \bideg v_i=(0,2),\quad\bideg u_i=(-1,2);\\
    d(1\otimes u_i)=\l_i\otimes 1,\quad d(v_i\otimes 1)=0,
  \end{align*}
\end{theorem}
Hence, in the case when $K$ is Cohen--Macaulay, the cohomology of
$U(K)$ can be calculated via the finite-dimensional differential algebra
$\k(K)/\!J\otimes\Lambda[u_1,\ldots,u_{m-n}]$ instead of infinite-dimensional
algebra $\k(K)\otimes\Lambda[u_1,\ldots,u_m]$ from Theorem~\ref{cohom2}.

\begin{example}
Let $K$ be the boundary complex of an $(m-1)$-dimensional
simplex. Then $\k(K)=\k[v_1,\ldots,v_m]/(v_1\cdots v_m)$. It easy to check
that only non-trivial cohomology classes in
$H\bigl[\k(K)\otimes\Lambda[u_1,\ldots,u_m],d\bigr]$ (see
Theorem~\ref{cohom2}) are represented by the cocycles $1$ and
$v_1v_2\cdots v_{m-1}u_m$ or their multiples. We have $\deg
(v_1v_2\cdots v_{m-1}u_m)=2m-1$, and Proposition~\ref{poin} shows that
$v_1v_2\cdots v_{m-1}u_m$ is the fundamental cohomological class of $\zk\cong
S^{2m-1}$ (see Example~\ref{zp}~1)).
\end{example}

\begin{example}
Let $K$ be a disjoint union of $m$ vertices. Then $U(K)$
is obtained by removing from $\C^m$ all codimension-two coordinate
subspaces $z_i=z_j=0$, $i,j=1,\ldots,m$ (see Example~\ref{uk}), and
$\k(K)=\k[v_1,\ldots,v_m]/I$, where $I$ is the ideal generated by all
monomials $v_iv_j$, $i\ne j$. It is easily deduced from
Theorem~\ref{cohom2} and Proposition~\ref{monom} that any cohomology class
of $H^{*}\bigl(U(K)\bigr)$ is represented by
a linear combination of monomial cocycles
$v_{i_1}u_{i_2}u_{i_3}\cdots u_{i_k}\subset
\k(K)\otimes\Lambda[u_1,\ldots,u_m]$ such that $k\ge2$, $i_p\ne i_q$ for
$p\ne q$. For each $k$ there $m\binom{m-1}{k-1}$ such monomials, and
there $\binom{m}{k}$ relations between them (each relation is obtained by
calculating the differential of $u_{i_1}\cdots u_{i_k}$).
Since $\deg(v_{i_1}u_{i_2}u_{i_3}\cdots u_{i_k})=k+1$, we have
\begin{align*}
  &\dim H^{0}\bigl(U(K)\bigr)=1,\quad
  H^{1}\bigl(U(K)\bigr)=H^{2}\bigl(U(K)\bigr)=0,\\
  &\dim H^{k+1}\bigl(U(K)\bigr)=\textstyle
  m\binom{m-1}{k-1}-\binom{m}{k},\;2\le k\le m,
\end{align*}
and the multiplication in the cohomology is trivial.

In particular, for $m=3$ we have 6 three-dimensional cohomology classes
$v_iu_j$, $i\ne j$, with 3 relations $v_iu_j=v_ju_i$, and 3
four-dimensional cohomology classes $v_1u_2u_3$, $v_2u_1u_3$, $v_3u_1u_2$
with one relation
$$
  v_1u_2u_3-v_2u_1u_3+v_3u_1u_2=0.
$$
Hence, $\dim H^{3}\bigl(U(K)\bigr)=3$, $\dim H^{4}\bigl(U(K)\bigr)=2$,
and the multiplication is trivial.
\end{example}

\begin{example}
Let $K$ be a boundary complex of an $m$-gon ($m\ge4$). Then, as it have been
mentioned above, the moment-angle complex
$\zk$ is a smooth manifold of dimension $m+2$, and $U(K)$
is homotopically equivalent to $\zk$. We have
$\k(K)=\k[v_1,\ldots,v_m]/I$, where $I$ is generated by monomials $v_iv_j$
such that $i\ne j\pm1$. (Here we use the agreement $v_{m+i}=v_i$ and
$v_{i-m}=v_i$.) The cohomology rings of these manifolds were calculated
in~\cite{BP2}. We have
$$
  \dim H^k\bigl(U(K)\bigr)=\left\{
  \begin{array}{l}
    1\quad\text{if }k=0\text{ or }m+2;\\[1mm]
    0\quad\text{if }k=1,2,m\text{ or }m+1;\\[1mm]
    (m-2)\binom{m-2}{k-2}-\binom{m-2}{k-1}-\binom{m-2}{k-3}
    \quad\text{if }3\le k\le m-1.
  \end{array}
  \right.
$$
For example, in the case $m=5$ there 5 generators of $H^3\bigl(U(K)\bigr)$
represented by the cocycles
$v_iu_{i+2}\in\k(K)\otimes\Lambda[u_1,\ldots,u_5]$, $i=1,\ldots,5$,
and 5 generators of $H^4\bigl(U(K)\bigr)$ represented by the cocycles
$v_ju_{j+2}u_{j+3}$, $j=1,\ldots,5$. As it follows from
Proposition~\ref{poin}, the product of cocycles $v_iu_{i+2}$ and
$v_ju_{j+2}u_{j+3}$ represents a non-trivial cohomology class in
$H^7\bigl(U(K)\bigr)$ (the fundamental cohomology class up to sign) if and
only if $\{i,i+2,j,j+2,j+3\}=\{1,2,3,4,5\}$.
Hence, for each cohomology class $[v_iu_{i+2}]$ there is a unique
(Poincar\'e dual) cohomology class $[v_ju_{j+2}u_{j+3}]$ such that
the product $[v_iu_{i+2}]\cdot[v_ju_{j+2}u_{j+3}]$ is
non-trivial.
\end{example}


\begin{thebibliography}{BR2}

\bibitem[Ar]{Ar}
V.\,I.~Arnold,
{\it The cohomology ring of the colored braid group} (Russian),
Mat. Zametki {\bf 5} (1969), 227--231;
English transl. in:
Math. Notes {\bf 5} (1969), 138--140.

\bibitem[Au]{Au}
M.~Audin,
{\it The Topology of Torus Actions on Symplectic Manifolds},
Progress in Mathematics~{\bf 93}, Birkh\"auser, Boston Basel Berlin, 1991.

\bibitem[Ba]{Ba}
V.\,V.~Batyrev,
{\it Quantum Cohomology Rings of Toric Manifolds},
Journ\'ees de G\'eom\'etrie Alg\'ebrique d'Orsay (Juillet 1992),
Ast\'erisque {\bf 218}, Soci\'ete Math\'ematique de France, Paris,
1993, pp. 9--34;
available at http://xxx.lanl.gov/abs/alg-geom/9310004.

\bibitem[Br]{Br}
E.~Brieskorn,
{\it Sur le groupes de tresses},
in: S\'eminare Bourbaki 1971/72, Lecture Notes in Math. {\bf 317},
Springer-Verlag, Berlin-New York, 1973, pp. 21--44.

\bibitem[BP1]{BP1}
V.\,M.~Bukhshtaber and T.\,E.~Panov,
{\it Algebraic topology of manifolds defined by simple polytopes}
(Russian), Uspekhi Mat. Nauk {\bf 53} (1998), no.~3, 195--196;
English transl. in:
Russian Math. Surveys {\bf 53} (1998), no.~3, 623--625.

\bibitem[BP2]{BP2}
V.\,M.~Buchstaber and T.\,E.~Panov,
{\it Torus actions and combinatorics of polytopes} (Russian),
Trudy Matematicheskogo Instituta im. Steklova {\bf 225} (1999), 96--131;
English transl. in: Proceedings of the Steklov Institute of
Mathematics~{\bf 225} (1999), 87--120;
available at http://xxx.lanl.gov/abs/math.AT/9909166.

\bibitem[BR1]{BR1}
V.\,M.~Bukhshtaber and N.~Ray,
{\it Toric manifolds and complex cobordisms} (Russian),
Uspekhi Mat. Nauk {\bf 53} (1998), no.~2, 139--140;
English transl. in:
Russian Math. Surveys {\bf 53} (1998), no.~2, 371--373.

\bibitem[BR2]{BR2}
V.\,M.~Buchstaber and N.~Ray,
{\it Tangential structures on toric manifolds, and connected sums of
polytopes}, preprint UMIST, Manchester, 1999.

\bibitem[Co]{Co}
D.\,A.~Cox,
{\it Recent developments in toric geometry},
in: Algebraic geometry (Proceedings of the Summer
Research Institute, Santa Cruz, CA, USA, July 9--29, 1995), J.~Kollar,
(ed.) et al. Providence, RI: American Mathematical Society. Proc. Symp.
Pure Math. {\bf 62} (pt.2), 389-436 (1997);
available at http://xxx.lanl.gov/abs/alg-geom/9606016.

\bibitem[Da]{Da}
V.~Danilov,
{\it The geometry of toric varieties},
(Russian), Uspekhi Mat. Nauk {\bf 33} (1978), no.~2, 85--134;
English transl. in:
Russian Math. Surveys {\bf 33} (1978), 97--154.

\bibitem[DJ]{DJ}
M.~Davis and T.~Januszkiewicz,
{\it Convex polytopes, Coxeter orbifolds and torus actions},
Duke Math. Journal {\bf 62}, (1991), no.~2, 417--451.

\bibitem[dCP]{dCP}
C. De Concini and C. Procesi,
{\it Wonderful models of subspace arrangements},
Selecta Mathematica, New Series {\bf 1} (1995), 459--494.

\bibitem[dL]{dL}
M. De Longueville,
{\it The ring structure on the cohomology of coordinate subspace
arrangements},
preprint, 1999;
available at http://www.math.tu-berlin.de/\~{}ziegler.

\bibitem[Fu]{Fu}
W.~Fulton,
{\it Introduction to Toric Varieties,}
Princeton Univ. Press, 1993.

\bibitem[GM]{GM}
M. Goresky and R. MacPherson,
{\it Stratified Morse Theory},
Springer-Verlag, Berlin-New York, 1988.

\bibitem[Ho]{Ho}
M. Hochster,
{\it Cohen--Macaulay rings, combinatorics, and simplicial complexes},
in: Ring Theory II (Proc. Second Oklahoma Conference),
B.\,R.~McDonald and R.~Morris, editors, Dekker, New York, 1977,
pp.~171--223.

\bibitem[Ma]{Ma}
S. Maclane,
{\it Homology},
Springer-Verlag, Berlin, 1963.

\bibitem[OS]{OS}
P.~Orlik and L.~Solomon,
{\it Combinatorics and Topology of Complements of Hyperplanes},
Invent. Math. {\bf 56} (1980), 167--189.

\bibitem[Pa1]{Pa1}
T.\,E.~Panov,
{\it Combinatorial formulae for the $\chi_y$-genus of a multioriented
quasitoric manifold} (Russian),
Uspekhi Mat. Nauk {\bf 54} (1999), no. 5, 169--170;
English translation in:
Russian Math. Surveys {\bf 54} (1999), no. 5.

\bibitem[Pa2]{Pa2}
T.\,E.~Panov,
{\it Hirzebruch genera of manifolds with torus action},
preprint, 1999;
available at http://xxx.lanl.gov/abs/math.AT/9910083.

\bibitem[Sm]{Sm}
L. Smith,
{\it Homological Algebra and the Eilenberg--Moore Spectral Sequence},
Transactions of American Math. Soc. {\bf 129} (1967), 58--93.

\bibitem[St]{St}
R. Stanley,
{\it Combinatorics and Commutative Algebra},
Progress in Math. {\bf 41}, Birkh\"auser, Boston, 1983.

\bibitem[Yu]{Yu}
S. Yuzvinsky,
{\it Small rational model of subspace complement},
preprint, 1999;
available at http://xxx.lanl.gov/abs/math.CO/9806143.


\end{thebibliography}
\end{document}